\documentclass{amsart}

\usepackage{amssymb}

\def\bc{\begin{center}}
\def\ec{\end{center}}
\def\be{\begin{equation}}
\def\ee{\end{equation}}
\def\ben{\begin{enumerate}}
\def\een{\end{enumerate}}
\def\bfg{\begin{figure}}
\def\efg{\end{figure}}
\def\bq{\begin{quote}}
\def\eq{\end{quote}}
\def\bd{\begin{description}}
\def\ed{\end{description}}
\def\this{i.\ e.\ } 

\def\h{\hbar}
\def\p{\partial}
\def\w{\wedge}

\def\dim{\operatorname{dim}}

\def\det{\operatorname{det}}

\newcommand{\CC}{{\Bbb C}}

\newcommand{\QQ}{{\Bbb Q}}
\newcommand{\EE}{{\Bbb E}}

\newcommand{\lan}{\langle}
\newcommand{\ran}{\rangle}

\newcommand{\gx}{\chi}

\renewcommand{\L}{\Lambda}
\renewcommand{\l}{\lambda}

\newcommand{\M}{\overline{\mathcal M}}
\newcommand{\ct}{\operatorname{ct}}
\newcommand{\ev}{\operatorname{ev}}
\renewcommand{\QQ}{\mathbf Q}
\renewcommand{\CC}{\mathbf C}
\renewcommand{\t}{\mathbf t}
\newcommand{\q}{\mathbf q}

\newcommand{\D}{\mathcal D}
\newcommand{\E}{\mathcal E}
\newcommand{\A}{\mathcal A}
\renewcommand{\a}{\alpha}
\renewcommand{\b}{\beta}
\renewcommand{\c}{\gamma}
\renewcommand{\d}{\delta}

\renewcommand{\H}{{\mathcal H}}
\newcommand{\F}{{\mathcal F}}
\newcommand{\G}{{\mathcal G}}
\newcommand{\T}{{\mathcal T}}
\newcommand{\C}{{\mathcal C}}
\renewcommand{\L}{{\mathcal L}}
\renewcommand{\O}{{\Omega}}
\newcommand{\1}{{\bf 1}}
\newcommand{\s}{{\mathbf s}}
\renewcommand{\S}{{\Sigma}}

\begin{document}

\title[Gromov -- Witten invariants and quantization]
{Gromov -- Witten invariants and quantization of quadratic hamiltonians}
 
\author{Alexander B. Givental} 
\address{UC Berkeley and Caltech} 
\thanks{Research partially supported by NSF Grant DMS-0072658} 

\date{April $12$, $2001$, revised on July $15$ and August $14$, 2001.}

\begin{abstract}
We describe a formalism based on 
quantization of quadratic hamiltonians and symplectic actions of loop groups 
which provides a
convenient home for most of known general results and conjectures about
Gromov-Witten invariants of compact symplectic manifolds and, more
generally, Frobenius structures at higher genus. We state several
results illustrating the formalism and its use. In particular, 
we establish Virasoro 
constraints for semisimple Frobenius structures and outline a proof
of the Virasoro conjecture for Gromov -- Witten invariants of complex 
projective spaces and other Fano toric manifolds. 
Details will be published elsewhere.   
\end{abstract}

\maketitle
 
{\large \bf 1. Gromov -- Witten invariants.}

\medskip
 
Let $X$ be a compact almost K\"ahler manifold of complex dimension $D$.
Denote by $X_{g,m,d}$ the moduli (orbi)space of degree $d$ stable holomorphic 
maps to $X$ of genus $g$ curves with $m$ marked points \cite{Kn1, BM}. 
The degree
$d$ takes values in the lattice $H_2(X)$. The moduli space is compact and
can be equipped \cite{BF, LT, S}  with a rational coefficient 
{\em virtual fundamental cycle} 
$[X_{g,m,d}]$ of complex dimension $m+(1-g)(D-3)+\int_{d} c_1(TX)$.

The {\em total descendent potential} of $X$ is defined as
\[ {\mathcal D}_X := \exp \sum \h^{g-1} {\mathcal F}^g_X ,\]
where ${\mathcal F}^g_X$ is the {\em genus $g$ descendent potential}
\[ \sum_{m,d} \frac{Q^d}{m!} \int_{[X_{g,m,d}]} 
 \wedge_{i=1}^{m}(\sum_{k=0}^{\infty} (\ev_i^*t_k) \psi_i^k ) .\]
Here $\psi_i^k$ are the powers of the $1$-st Chern class of the universal 
cotangent line bundle over $X_{g,m,d}$ corresponding to the $i$-th marked 
point, $\ev_i^*t_k$ are pull-backs by the evaluation map 
$\ev_i: X_{g,m,d}\to X$ at the $i$-th marked point of the cohomology classes
$t_0,t_1,... \in H^*(X,\QQ)$, and $Q^d$ is the representative of $d$ in 
the semigroup ring of the semigroup of degrees of holomorphic curves in $X$.  
The {\em genus $g$ Gromov -- Witten potential} of $X$ is defined as the
restriction $F_X^g(t):= \F^g_X |_{t_0=t, \ t_1=t_2=...=0} $.

The genus $g$ descendent potentials are considered as formal functions
on the super-space of vector Laurent polynomials $\t (z) =t_0+t_1z+t_2z^2+...$ 
in one indeterminate $z$ with coefficients in $H:=H^*(X; \QQ [[Q]])$,
the cohomology space of $X$ over the {\em Novikov ring}. The latter is
a suitable power series completion of the semigroup ring of degrees.          

\medskip

{\bf 1.1.} {\em Remark.} We will assume further on that $H$ has no odd part
and leave the super-space generalization of the material of sections 
$2$ -- $5$ to the reader. It is not clear at the moment if 
the content of sections $6$ -- $10$ admits such a generalization. 
 
\medskip

{\bf 1.2.} {\em Example.} 
When $X = point$, the moduli spaces of stable
maps coincide with the Deligne -- Mumford compactifications $\M_{g,m}$ 
of moduli spaces of genus $g$ Riemann surfaces with $m$ marked points.
According to Witten's conjecture \cite{W} proved by Kontsevich \cite{Kn},
the total descendent potential in this case coincides with the tau-function
of the KdV-hierarchy satisfying the string equation \cite{KS}.
We will denote ${\mathcal D}_{point}$ by $\tau (\h ; \t )$ and call it the
{\em Witten -- Kontsevich tau-function}.

\medskip

Taylor coefficients of the Gromov -- Witten and descendent potentials 
 depend only on the deformation class of the symplectic 
structure and are often referred to respectively as Gromov -- Witten 
invariants and their gravitational descendents.  

\bigskip     

{\large \bf 2. Quadratic hamiltonians and quantization.}

\medskip

Let $H, ( \cdot, \cdot )$ be an $N$-dimensional vector space  
equipped with a non-degenerate symmetric bilinear form. Let $\H$ be the space
of Laurent polynomials in one indeterminate $z$ with vector coefficients 
from $H$. We introduce a symplectic bilinear form in $\H$ by
\[ \O (f,g) = \frac{1}{2\pi i}\oint ( f (-z), g(z) ) dz = - \O (g,f) .\]

Let $\H=\H_{+}\oplus \H_{-}$ correspond to the decomposition 
$f(z,z^{-1})=f_{+}(z)+f_{-}(1/z)/z$ of the Laurent polynomials into 
polynomial and polar parts. The subspaces $\H_{\pm}$ are Lagrangian.
In particular the projection $\H\to \H_{+}$ along $\H_{-}$ defines a 
polarization of $(\H, \O)$. We quantize infinitesimal symplectic 
transformations $L$ on $\H$ to
order $\leq 2$ linear differential operators $\hat{L}$. In a 
Darboux coordinate system $\{ p_{\a}, q_{\b} \}$ compatible with our 
decomposition $\H=\H_{+}\oplus \H_{-}$ we have
\[  (p_{\a}p_{\b})\hat{\ }=\h \p_{q_{\a}}\p_{q_{\b}},\ 
(p_{\a}q_{\b})\hat{\ }= q_{\b} \p_{q_{\a}},\ 
(q_{\a}q_{\b})\hat{\ }=q_{\a}q_{\b}/\h .\]
Note that $ [\hat{F},\hat{G}] =\{ F,G \} \hat{\ }
+{\mathcal C} (F,G)$
where the cocycle ${\mathcal C}$ satisfies
\[ {\mathcal C}(p_{\a}^2, q_{\a}^2)=2,\ \ 
{\mathcal C} (p_{\a}p_{\b},q_{\a}q_{\b})=1 \ \text{for}\ \a\neq \b,\]  
and ${\mathcal C} =0$ for any other pairs of quadratic Darboux monomials. 

\medskip

{\bf 2.1.} {\em Remark: Fock spaces.} 
The differential operators $\hat{L}$ act on functions of $\h$ and  
${\bf q} =  q_0+q_1z+q_2z^2+...$ where $q_0,q_1,q_2,... \in H$. 
We will often refer to such functions as elements of the {\em Fock space}.
However we will be concerned with various functions which belong
to different formal series completions of the polynomial Fock space,
and we will not describe the completions explicitly in this paper.

\medskip

{\bf 2.2.} {\em Dilaton shift.} The total descendent potential 
${\mathcal D}_X$ will be considered as a vector in the Fock 
space corresponding to the space $H=H^*(X; \QQ [[Q]])$ equipped 
with the Poincar\'e intersection form $(\cdot ,\cdot )$. 
More precisely, we introduce the {\em dilaton shift} $\q (z) = \t (z)-z$
which identifies the indeterminates $t_0,t_1,...$ in the descendent 
potential with coordinates $q_0,q_1,... $ in $\H_{+}$. 
Slightly abusing notation we will always identify functions of $\t $ 
with the functions of $\q \in \H_{+}$ obtained from them by the dilaton shift.
In particular, the genus $g$ potentials ${\mathcal F}^g_X$ are considered
as formal functions of $\t$ near the origin $(0,0,0,...)$ or, equivalently,
as formal functions of $\q $ near $(0,-1,0,...)$. The same convention
applies to the Witten -- Kontsevich tau-function $\tau = {\mathcal D}_{point}$.

\bigskip          

{\large \bf 3. The Witten -- Kontsevich tau-function.}

\medskip
 
Denote $D=z(d/dz)z$ and put $L_m = z^{-1/2} D^{m+1} z^{-1/2}$ so that 
\[ L_{-1}= 1/z,\    L_{0}=z\frac{d}{dz} +1/2,\ 
L_1=z^3\frac{d^2}{dz^2}+3z^2\frac{d}{dz}+\frac{3}{4}z,\]
\[ L_2=z^5\frac{d^3}{dz^3}+\frac{15}{2}z^4\frac{d^2}{dz^2}+\frac{45}{4}z^3
\frac{d}{dz}+\frac{15}{8}z^2, \ ... \]
It is easy to check that $\O (Df,g)=\O (f, Dg)$ and 
that $\O (z^{-1/2}f,g)=-i\O(f,z^{-1/2}g)$ (whatever it means). This implies 
$\O(L_m f,g)=-\O (f,L_mg)$ and shows that the operators $L_m$
are infinitesimal symplectic transformations on $\H$. On the other hand, 
$D$ is conjugate to $z^2d/dz=-d/d w$ where $w=1/z$, and hence $L_m$ commute 
as $-wd^{m+1}/dw^{m+1}$ and therefore (via the Fourier transform)
--- as the vector fields $-x^{m+1}d/dx$ on the line. Thus the
Poisson brackets satisfy $\{L_m,L_n\}=(m-n)L_{m+n}$, and we have a 
representation of the Lie algebra of vector fields on the line to the Lie 
algebra of quadratic hamiltonians on $\H$. 

In the case of $1$-dimensional $H$ with the standard inner product, using 
the Darboux coordinate system $f=...+p_1/z^2-p_0/z+q_0+q_1z+...$ on $\H$,
we get (here $\p_k=\p/\p q_k$):
\[ \begin{array}{l}
 \hat{L}_{-1}=q_0^2/2\h+ \sum_{m\geq 0} q_{m+1}\p_m \\
 \hat{L}_0=  \sum_{m\geq 0} (m+1/2)q_m\p_m \\
 \hat{L}_1=\h \p_0^2/8+ \sum_{m\geq 0}(m+1/2)(m+3/2)q_m\p_{m+1} \\
 \hat{L}_2=3\h \p_0\p_1/4+\sum_{m\geq 0}(m+1/2)(m+3/2)(m+5/2)q_m\p_{m+2} .
\end{array} \]
We have $[\hat{L}_m,\hat{L}_n]=(m-n)\hat{L}_{m+n}$ unless $m,n=\pm 1$,
in which case: $[\hat{L}_1,\hat{L}_{-1}]=2[\hat{L}_0+1/16]$.
Thus the operators $\hat{L}_m+\d_{m,0}/16$ form a representation of the 
Lie algebra of vector fields on the line. 

Due to \cite{KS}, the following result is a reformulation of the
Kontsevich theorem \cite{Kn} confirming 
the Witten conjecture \cite{W}. 

\medskip

{\bf 3.1. Theorem.} {\em The Witten -- Kontsevich tau-function is annihilated
by the operators $\hat{L}_m+\d_{m,0}/16$, $m=-1,0,1,2,...$, and is completely
characterized by this property (up to a scalar factor).}

\bigskip
 
{\large \bf 4. Hodge integrals.}

\medskip
 
Let $\EE $ denote the {\em Hodge bundle} over the moduli space $X_{g,m,d}$.
By definition the fiber of $\EE$ over the point represented by a stable map 
$\S \to X$ is the complex space of dimension $g$ dual to 
$H^1(\S , {\mathcal O}_{\S})$. 
In fact ${\EE}$ is the pull-back of by the 
{\em contraction map}  $\ct: X_{g,m,d}\to \M_{g,m}$ of the
Hodge bundle over the Deligne -- Mumford space. 
It is known \cite{Mu} that even components $ch_{2k}(\EE)$ 
of the Chern character vanish. 
We define the {\em total Hodge potential} of $X$ as an extension of 
the total descendent potential depending on the sequence 
$\s = (s_1,s_2,...)$ of new variables and incorporating intersection indices 
with characteristic classes of the Hodge bundles (as in \cite{FP}):
\[ \E_X (\h ; \t ; \s ):= \exp \sum_g \h^{g-1} 
\sum_{m,d} \frac{Q^d}{m!}  \int_{[X_{g,m,d}]} 
e^{ \sum_{k=1}^{\infty} s_k\ ch_{2k-1}(\EE)}
\wedge_{i=1}^m (\sum_k (\ev_i^*t_k) \psi_i^k ) . \]  

We consider $\E_X$, subject to the dilaton shift, as a family of elements 
in the Fock space depending on the parameters $\s$.

On the other hand, it is obvious that multiplication by 
$z^{2k-1}$ defines
an infinitesimal symplectic transformation on  
$(\H, \O )$, and we denote by $(z^{2k-1})\hat{\ }$ the corresponding 
quantization.

\medskip

{\bf 4.1. Theorem.} 
\[ \E_X \  = \  \exp [ \sum_{k=1}^{\infty} 
\frac{B_{2k}}{(2k)!}\ s_k\ (z^{2k-1})\hat{\ } ] \ \ {\D}_X  \]   
{\em where $B_{2k}$ are Bernoulli numbers: $x/(1-e^{-x})=
1+x/2+\sum_{k=1}^{\infty} B_{2k} x^{2k}/(2k)!$.}

\medskip

{\bf 4.2.} {\em Remark.} The theorem is a reformulation of a result
by Faber -- Pandharipande  
\cite{FP} which adjusts to arbitrary target spaces the famous Mumford's 
Riemann -- Roch -- Grothendieck formula \cite{Mu} expressing 
$ch (\EE )$ via $\psi_i$.    

\bigskip

{\large \bf 5. Gravitational ancestors.} 

\medskip

Consider the composition $\pi: X_{g,m+l,d}\to \M_{g,m+l}$ $\to \M_{g,m}$
of the contraction map with the operation of forgetting all marked points
except the first $m$. Let $\bar{\psi}_i:=\pi^*(\psi_i)$ denote pull-backs
of the classes $\psi_i$, $i=1,...,m$, from $\M_{g,m}$. We introduce
the {\em total ancestor potential} 
\[ \A_t (\h; \t):= \exp \sum_{g=0}^{\infty} \h^{g-1}\bar{\F}^g_t ,\]
where the {\em genus $g$ ancestor potential $\bar{\F}^g_t$} is defined by
\[ \bar{\F}^g_t:=\sum_{m,l,d}\frac{Q^d}{m! l!}
\int_{[X_{g,m+l,d}]}\wedge_{i=1}^m (\sum_k (\ev^*_it_k) \bar{\psi}_i^k) \ 
\wedge_{i=m+1}^{m+l} \ev_i^*t .\]
By definition the sum does not contain the terms with $(g,m)=(0,0),(0,1),(0,2)$
and $(1,0)$.
We treat $\A_t$ subject to the dilaton shift $\q(z) = \t (z)-z \in \H_{+}$ 
as an element in the Fock space depending on the parameter $t\in H$.
   
Recall that $F^1(t)$ denotes the genus $1$ Gromov -- Witten potential, \this 
$\F^1_X(\t) $ at $t_0=t, t_1=t_2=...=0$.

\medskip
 
Let us introduce the operator $S_t$ on the Laurent 
$1/z$-series completion of the space $H$ defined by
\[ (a , S_t b):= \lan a, \frac{b}{z-\psi}\ran := 
(a,b)+\sum_{k=0}^{\infty} \lan a, b \psi^k\ran \ z^{-1-k}.\]
The correlator notation here refers to the following {\em $2$-point
gravitational descendent} in genus $0$: 
\[ \lan a \psi^k, b \psi^l\ran := \sum_{m,d}\frac{Q^d}{m!}\int_{[X_{0,m+2,d}]}
(\ev_1^*a) \psi_1^k \ (\wedge_{i=2}^{m+1}\ev_i^*t\ )\ (\ev_{m+2}^*b)\ \psi_{m+2}^l .\]
It is one of the basic facts of quantum cohomology theory (see the next
section) that $S_t^*(-1/z)S_t(1/z)=\1 $. (The asterisk here means 
transposition with respect to the inner product $(\cdot ,\cdot )$ on $H$.)
In other words, on a suitable completion of $(\H, \O)$, the operator $S_t$ defines 
a symplectic transformation depending on the parameter $t\in H$.   
We put $\hat{S}_t := \exp (\ln S_t)\hat{\ }$. 

\medskip

{\bf 5.1. Theorem.} 
\[ \D_X = e^{F^1(t)} \hat{S}^{-1}_t \ \A_t. \]  

\medskip

{\bf 5.2.} {\em Remark.} The theorem is a reformulation of results by Kontsevich -- Manin \cite{KM} computing gravitational ancestors in terms of descendents and vice
versa. 

Moreover, using these results, Getzler \cite{Ge} 
proves the {\em $3g-2$-jet conjecture} of Eguchi
-- Xiong and Dubrovin about genus $g$ descendent potential $\F^g_X$ 
(we are not going to formulate it here --- see \cite{EX} for details)
by showing that it is equivalent to the following property of the genus 
$g$ ancestor potentials:
\[ \frac{\p^m }{\p t^{\a_1}_{k_1+1} ... \p t^{\a_m}_{k_m+1}} \ 
( \bar{\F}^g_t )|_{t_0=0} =0 \ \ \text{if} \ \ k_1+...+k_m > 3g-3.\]
The latter property is obvious:
$\bar{\psi}_1^{k_1+1}...\bar{\psi}_m^{k_m+1}=0$ 
since  $\dim \M_{g,m}=3g-3+m$. 

In particular, $\bar{\F}^0_t|_{t_0=0}=0$.  
One can use this property 
in order to extract the genus $0$ descendent potential $\F^0$ as follows.

\medskip

{\bf 5.3. Proposition.} {\em Quantized symplectic operators of the form 
$S(1/z)=\1+S_1/z+S_2/z^2+...$ act on elements $\G$ of the Fock space as   
\[ (\hat{S}^{-1} \G ) (\q ) = e^{ W(\q,\q) /2\h } \G ( [S \q ]_{+} ) ,\]
where $[ S \q ]_{+} $ is the power series truncation of $S(z^{-1}) \q (z)$,
and the quadratic form $W =\sum ( W_{kl} q_k, q_l)$ is defined by
\[ \sum_{k,l\geq 0} \frac{W_{kl}}{w^k z^l} := \frac{S^*(w^{-1})S(z^{-1}) -\1 }
{w^{-1}+z^{-1}} .\]}

\medskip

In particular, $\F^0_X (\q) = W_t(\q,\q)/2$ if the parameter $t$ in $S_t$ is
set to make the ancestor variable $t_0=0$.  

\medskip  

{\bf 5.4. Corollary.} {\em The genus $0$ descendent potential equals 
\[ \F^0_X (\t ) = \frac{1}{2} \lan \t (\psi)-\psi, \t(\psi)-\psi \ran\ |_{t=t(\t )} \]
where $t(\t)$ is the critical point of the function 
$\lan \1, \t(\psi)-\psi \ran  $ of $t\in H$ depending on the parameters 
$\t = (t_0,t_1,...)$.}

\medskip

{5.5.} {\em Remarks.} The corollary is the famous reconstruction results 
for genus $0$  gravitational
descendents due to Dubrovin \cite{D} 
and Dijkgraaf -- Witten \cite{DW}.

The ``mysterious'' requirement that $t=t(\t)$ is the critical point of 
$\lan \1, \t(\psi)-\psi\ran $ arises here simply to set the argument 
$ [S_t \q]_{+}$ in the ancestor potential equal $0$ at $z=0$.

\medskip

We would like to mention the following heuristic observation (one can make
on the basis of $5.3$, the dual proposition $7.3$ and the theorems $4.1$, 
$5.1$ and underlying geometry) about the nature of out formulas. 
It appears that their four basic ingredients --- symplectic transformations, 
quantization, dilaton shift and the central charge
--- are governed by the four missing Deligne -- Mumford spaces with
$(g,m)\ =\ (0,1),\ (0,2),\ (0,0)$ and $(1,0)$.
    
\bigskip
    
{\large \bf 6. Frobenius structures.} 

\medskip

The operator series $S_t(1/z)$ introduced in the
previous section is known to have the following properties \cite{Gi2,Gi1,Ma}. 

The operator-valued $1$-form $A(t):=z\ (d_t S_t(1/z))\  S_t^{-1}(1/z)$
does not depend on $z$ and thus defines a linear pencil of connections
$\nabla_z := d  - z^{-1}A(t) \w $ on the tangent bundle $TH$ flat for all 
values of the parameter $z^{-1}$. The flatness condition thus reads
$A\w A=0$ and $dA=0$. In coordinates $t=\sum t^{\a} \phi_{\a} $ the
first condition means commutativity $A_{\a}A_{\b}=A_{\b}A_{\a}$ of the
components of $A=\sum A_{\a}(t) dt^{\a}$. The natural correspondence
$\p_{\a}\mapsto A_{\a}(t)=i_{\p_{\a}}A $ defines commutative associative
multiplications $\bullet_t$ on the tangent spaces $T_t H$ called the 
{\em quantum cup-product}. Its structure constants 
$( \phi_{\a} \bullet_t \phi_{\b}, \phi_{\c})$ actually coincide with the
third directional derivatives $\p_{\a}\p_{\b}\p_{\c} F^0_X(t)$ of the 
genus $0$ Gromov -- Witten potential $F^0_X$. In particular, this explains
why $dA=0$ and also shows that $A_{\a}^*=A_{\a}$ so that the quantum 
cup-product is {\em Frobenius}: $( a\bullet b,c)=(a,b\bullet c)$.
The quantum cup-product on $TH$ is invariant under the translations in
the direction of the unit element $\1 \in H=H^*(X;\QQ [[Q]])$ and 
$\1 \bullet_t =\operatorname{id}$. 

The picture just described has been axiomatized \cite{D} under the name 
{\em Frobenius structure}. We refer to \cite{D, Ma, Gi2, Gi1, Gi0} for 
discussions of the following definition and derivations of the properties  
reviewed below.         

\medskip

{\bf 6.1.} {\em Definition.} A Frobenius structure on a
manifold $H$ consists of:

\noindent (i) a flat pseudo-Riemannian metric $(\cdot ,\cdot)$,
    
\noindent (ii) a function $F$ whose $3$-rd covariant derivatives 
$F_{abc}$ are structure constants $(a\bullet_t b,c)$ of a Frobenius algebra
structure, i.e. associative commutative multiplication
$\bullet_t $ satisfying $(a\bullet_t b, c)=(a,b\bullet_t c)$, 
on the tangent spaces $T_tH$ which depends smoothly on $t$; 

\noindent (iii) the vector field of unities ${\1}$ of the $\bullet_t$-product
which has to be covariantly constant and preserve the multiplication 
and the metric.  

\medskip

In the flat affine coordinates $t=\sum t^{\a} \phi_{\a}$ of the metric 
$(\cdot, \cdot)$, consider the deformation 
$\nabla_z:=d-z^{-1}\sum (\phi_{\a}\bullet_t) \ d\t_{\a}\w $ of the Levi-Civita
connection. It follows from the definition that $\nabla_z^2=0$ for all $z$.
Furthermore, one can construct a fundamental solution to the system
$\nabla_z S = 0$ in the form of a power series $S(1/z)=\1+S_1/z+S_2/z^2+...$
satisfying $S^*(-1/z)S(1/z)=\1 $. Such $S$ is unique up to right multiplication
by a constant operator series $C(1/z)=\1+C_1/z+C_2/z^2+...$ satisfying
$C^*(-1/z)C(1/z)=\1$. We will call a Frobenius manifold equipped with
a choice of the solution $S$ {\em calibrated}. 

\medskip

{\bf 6.3.} {\em The Euler field}  on a Frobenius
manifold is a vector field $E$ which in a flat affine coordinate system 
on $H$ has the form of a linear inhomogeneous vector field and such that 
$\bullet , \1 $ and $(\cdot ,\cdot )$ are eigenvectors
of the Lie derivative along $E$ with the eigenvalues $0$, $-1$ and $2-D$
respectively. The Frobenius manifold equipped with an Euler vector field 
is called {\em conformal} of conformal dimension $D$. We will require
that calibrations $S$ of conformal Frobenius structures are 
homogeneous in the sense that the bilinear form $(a, S(1/z) b)$ on $TH$ 
is an eigenvector of $z\p_{z}+E$ with the eigenvalue $2-D$. This
reduces the choice of $S$ to finitely many constants, and we refer to 
\cite{D1} for precise description of the ambiguity.
 
Following \cite{D1}, we will assume the linear part of the Euler field 
semisimple. 

\medskip

{\bf 6.4.} {\em Examples.} (a) According to K. Saito \cite{Sa} orbit
spaces of finite irreducible Coxeter groups carry canonical conformal
Frobenius structure. 

(b) Miniversal deformations of isolated critical points of holomorphic 
functions can be naturally equipped with conformal Frobenius structures
\cite{Sa1} (see also \cite{Ma,Gi3,Gi0}). 

(c) Genus $0$ Gromov -- Witten invariants of a compact almost K\"ahler manifold
$X$ define on $H=H^*(X;\QQ[[Q]])$ a formal structure of a calibrated 
conformal Frobenius manifold of conformal dimension $D=\dim_{\CC} X$. 
In the coordinate system $t=\sum t^{\a}\phi_{\a}$
corresponding to a graded basis $\{ \phi_{\a} \}$ in $H^*(X,\QQ)$, the
Euler field takes on the form
\[ E = \sum (1-\deg \phi_{a}/2) t_{\a} \p_{\a} + \rho \] 
where the constant part $\rho \in H^*(X)$ is the $1$-st Chern class of
the tangent bundle $TX$.

(d) In the section $9$ we will deal with equivariant generalization 
\cite{Gi2, S1} of Gromov -- Witten theory in the case when $X$ is equipped with
a hamiltonian action of a compact Lie group $G$. In this case equivariant
genus $0$ Gromov -- Witten invariants define on
equivariant cohomology space $H=H^*_G(X;\QQ [[Q]])$ the structure of a
calibrated Frobenius manifold. It is not conformal though since the Euler
field defined by cohomology grading over $\QQ [[Q]]$ is not a derivation
over the coefficient ring $H^*_G(pt) = H^*(BG)$ of equivariant cohomology
theory.   
   
\medskip

{\bf 6.5.} {\em Remark.} We will apply our formalism of quantized quadratic
hamiltonians to the general problem of equipping Frobenius manifolds with all attributes
of higher genus Gromov -- Witten theory. The function $F$ from the definition
$6.1$ is taken of course on the role of the genus $0$ potential $F^0$. The
calibration $S$ is then used to define the $1$-point descendent 
$\lan a,  b/(z-\psi )\ran := (a,S(1/z)b)$. Then Corollary $5.4$ yields
Dubrovin's construction \cite{D} of the genus $0$ descendent potential $\F^0$.

Below we give a construction of the total descendent potentials $\D $ 
in the case when the Frobenius manifold is {\em semisimple} \this when
Frobenius algebras $(T_tH, \bullet_t)$ are semisimple at generic $t\in H$.
All Frobenius manifolds of Examples $6.4 a,b$ are semisimple as well
as those in $d$ when $G$ is a torus acting on $X$ with isolated fixed points. 
Among examples $6.4c$ flag manifolds, toric Fano manifolds and 
probably many other Fano manifolds yield semisimple Frobenius structures. 

\medskip

{\bf 6.6.} {\em Canonical coordinates.} Let $t=\sum t^{\a}\phi_{\a}$ be
a flat coordinate system on the Frobenius manifold $H$. 
Let $u\in H$ be a semisimple point so that the matrices 
$A_{\a}=\phi_{\a}\bullet_{u}$ are simultaneously diagonalizable.
Let $\Psi $ be the transition matrix from the basis 
$\{ \phi_{\a} \}$ in $T_uH$ to the basis of common eigenvectors of
$A_{\a}$ normalized to the unit lengths. Then  $\Psi^{-1} A^1 \Psi $
is a diagonal matrix of {\em closed} $1$-forms and has locally the form
$dU =\operatorname{diag}(du^1,...,du^N)$. The diagonal entries of the potential
matrix  $U=\operatorname{diag}(u^1,...,u^N)$ 
form a local coordinate system on $H$ called 
{\em canonical} \cite{D}.
Canonical coordinates are defined uniquely up to signs, permutations and
additive constants which are set to zero in the conformal case by the
requirements $E u^i=u_i$.   

\medskip

{\bf 6.7. Proposition} \cite{Gi0,Gi1,D}.  
{\em (a) The equation $\nabla_{z} S =0$ in a neighborhood of a semisimple 
point $u$ has a fundamental solution in the form: $ \Psi_u R_u(z) e^{U/z} $, where
$R_u(z)=\1+R_1z+R_2z^2+...$ is a formal matrix power series
satisfying $R_u^*(-z)R_u(z)=\1$.

(b) The series $R_u(z)$ satisfying (a) 
is unique up to right multiplication by diagonal matrices 
$\exp (a_1z+a_2z^3+a_3z^5+ ...)$  where
$a_{k}=\operatorname{diag}(a_{k}^1,...,a_{k}^N)$ are constant.

(c) In the case of conformal Frobenius structures the
series $R_u(z)$ satisfying (a) is uniquely determined by
the homogeneity condition $(z\p_{z}+\sum u^i\p_{u^i}) R_u(z) =0$.}

\medskip
  
Let $\H $ be the space of Laurent polynomials in $z$
with coefficients in the tangent space $T_uH$ to the Frobenius manifold
at a semisimple point $u$. Let $\H=\H_{+}\oplus \H_{-}$ be the
polarization of $\H$ described in section $2$. For $\q (z) \in \H_{+}$
we write $\Psi_u^{-1}\q(z) = (\q^1(z),...,\q^N(z))$ and introduce the
direct product $ \T = \tau (\h ; \q^1 )...\tau (\h; \q^N)$ of $N$ copies
of the Witten -- Kontsevich tau-function as an element of the
Fock space of functions on $\H_{+}$. Next, the series $ R_u(z) $ defines
a symplectic transformation on $\H$, and we put 
$\hat{R}_u =\exp (\ln R_u)\hat{\ }$. Slightly abusing notation we denote
$\hat{\Psi}_u$ the operator $\G (\Psi^{-1} \q ) \mapsto \G (\q) $ identifying
the Fock space with its coordinate version. Finally, 
$\sum_{i=1}^N R_1^{ii}du^i $ where $R_1^{ii}$ are the diagonal entries
of the matrix $R_1$ in the series $R_u(z)=\1+R_1/z+...$ is known to be closed 
as a $1$-form on the Frobenius manifold. We introduce the function 
$C(u)=\frac{1}{2}\int^u \sum R_1^{ii} du^i$ of $u$ defined up to an additive constant.
 
\medskip

{\bf 6.8.} {\em Definition.} We define the {\em total descendent potential} of
a semisimple Frobenius manifold by the formula
\[ \D (\h; \t)\ :=\ e^{C(u)}\ \hat{S}_u^{-1}\ \hat{\Psi}\ \hat{R}_u \ 
e^{(U/z)\hat{ }}\ \T .\] 
We introduce the {\em total ancestor potential} of a semisimple
Frobenius manifold:
\[ e^{F^1(u)}\ \A_u (\h; \t)\ :=\ e^{C(u)}\ \hat{\Psi} \ \hat{R}_u \ 
e^{(U/z)\hat{ }}\ \T .\]
Since $\D $ and $\A_u$ are related as in Theorem $5.1$, these definitions
automatically agree with the reconstruction formula $5.4$ 
for $\F^0$. We keep the convention about the dilaton
shift $2.2$ in these definitions.  

\medskip

{\bf 6.9.} {\em Remarks.} (a) Both potentials depend on a choice of the 
asymptotical series $R_u(z)$ (which is unambiguous in the conformal case)
and are defined up to a constant factor {\em independent on $u$}.             
In section $9$ we will specify the choice of $R_u(z)$ corresponding to
equivariant Gromov -- Witten theory.  

(b) The factor $\exp (U/z)\hat{ }$ is redundant since 
the string operator $(1/z)\hat{ }=\hat{L}_{-1}$ annihilates the (product of) 
Witten -- Kontsevich tau-functions as it follows from $3.1$ and is included
only for future convenience. 

(c) Quantizations of symplectic transformations $S_u(1/z)$ and $R_u(z)$
require different power series completion of the space $\H_{+}$ so that
their composition, strictly speaking, is not defined. Nevertheless the
formula for $\D$ makes sense (in particular --- due to some properties 
of the function $\T$) at least in the formal neighborhood of the critical 
point locus $u=t(\t)$ described in $5.4$. As a result of this subtlety,
the potentials of abstract Frobenius structures thus defined do not always
extend to non-semisimple values of $t_0=t({\bf 0})$. We are not going to
stress this subtlety in the rest of the text.

\bigskip

{\large \bf 7. Properties of the total descendent potential.}

\medskip

{\bf 7.1. Theorem.} {\em The total descendent potential $\D $ of a semisimple
Frobenius manifold defined in $6.8$ does not depend on the choice of 
a semisimple point $u$.}

\medskip

This follows from the property of $S_u(1/z)$ and $\Psi_u R_u(z) \exp (U/z)$
to satisfy the same differential equation $\nabla_{z} S =0$ with coefficients
rational in $z$. As a result of this, derivatives of $\D$ in the directions 
of the parameter $u$ vanish as if $S_u^{-1}$ and $\Psi_u R_u \exp (U/z)$
were inverse to each other. The factor $e^{C(u)}$ however is needed
in order to offset the effect of the cocycle $\C $ under quantization of this 
computation.

\medskip

{\bf 7.2. Corollary.} {\em (a) The genus $1$ Gromov -- Witten potential $F^1(t)$ of
a semisimple Frobenius manifold is given by the formula
\[ F^1(t)= C(u) - \frac{1}{48} \ln \det (\p_{u^i},\p_{u^j}) =
\frac{1}{2}\int^u \sum  R_1^{ii}du^i + \frac{1}{48} \sum \ln \Delta_i  ,\]
where $\Delta_i^{-1}(u)=(\p_{u^i},\p_{u^i})$ are inner squares of the canonical
idempotents $\p_{u^i}$ of the semisimple Frobenius multiplication 
$\bullet_u$.

(b) The genus $1$ descendent potential equals
\[ \F^1_X(\t ) = F^1(t(\t))\ +\ \frac{1}{24} \ln \det 
[ \p_{t^{\a}_0}\p_{t^{\b}_0}\p_{t^{0}_0}\F_X^{0} (\t) ] ,\]
where the partial derivatives are taken with respect to
coordinates of $t^0=\sum t_0^{\a}\phi_{\a}$.
}
  
\medskip

The part (a) coincides with the conjectural formula for $F^1_X$ suggested in
\cite{Gi1} and proved in \cite{DZ} in the conformal case.
A similarly looking formula can be derived for the 
genus $1$ descendent potential. It is shown in \cite{Gi0}, Example $8$, 
how to reconcile it with the well-known
formula in part (b) which is due to
Dijgraaf -- Witten \cite{DW} and does not require
semisimplicity hypotheses.   

\medskip

{\bf 7.3. Proposition.} {\em Quantized symplectic operators of the form
$R(z)=1+R_1z+R_2z^2+...$ act on elements $\G $ of the Fock space as
\[ (\hat{R} \G) (\q) = [ e^{\h V(\p_{\q},\p_{\q})/2}\ \G ] (R^{-1} \q ) ,\]
where $R^{-1}\q$ is defined as the product of $z$-series $R^{-1}(z)\q(z)$, 
and the quadratic form $V=\sum (p_k, V_{kl} p_l)$ is defined by
\[ \sum_{k,l\geq 0} (-1)^{k+l}V_{kl} w^kz^l = \frac{R^*(w)R(z)-\1}{w+z} .\]}  
 
\medskip

{\bf 7.4. Corollary.} {\em The total descendent potential $\D$ of 
a semisimple Frobenius manifold satisfies the $3g-2$-jet conjecture
of Eguchi -- Xiong (see Remark $5.2$).}

\medskip 

This time the property of genus $g$ ancestor potentials
\[ \frac{\p^m}{\p_{t^{\a_1}_{k_1+1}}...\p_{t^{\a_m}_{k_m+1}}}\ 
(\bar{\F}^g_u )|_{t_0=0} \ =\ 0 \ \ \text{if} \ \ k_1+...+k_m>3g-3 \]
follows from the dimensional properties of the product $\T$ of
Witten-Kontsevich tau-functions and from the structure (``upper-triangular''
in some sense) of the operators $\hat{R}$ described in the proposition. 

Propositions $5.3$ and $7.3$ together imply:

\medskip 

{\bf 7.5. Corollary.} {\em The definition $6.8$ of the descendent 
potential of a semisimple Frobenius manifold agrees with the construction 
of genus $>1$ (descendent) potentials $F^g(t)$ (respectively $\F^g(\t)$) 
suggested in \cite{Gi0}.}

\medskip

{\bf 7.6. Conjecture.} {\em $\D_X=\D$: The total descendent potential $\D$ 
corresponding to the calibrated conformal Frobenius structure on 
the cohomology space $H=H^*(X;\QQ [[Q]])$ defined by genus $0$ Gromov --
Witten invariants of a compact almost K\"ahler manifold $X$ with
generically semisimple quantum cup-product coincides with the total 
descendent potential $\D_X$ defined in the Gromov -- Witten theory of $X$.}   

\medskip

This conjecture is equivalent to Conjecture $2$ in \cite{Gi0}. 
In sections $9, 10$ we outline a proof of this conjecture for 
complex projective spaces and other toric Fano manifolds. 

\medskip

The Virasoro operators $\hat{L}_m+N\d_{m,0}/16$ 
from the section $3$ annihilate 
the product $\T$ of $N$ Witten -- Kontsevich tau-functions. Let us 
put formally 
\[ \L_m := \hat{S}_u^{-1} \hat{\Psi}_u \hat{R}_u \ \hat{L}_m\ 
\hat{R}_u^{-1}\hat{\Psi}_u^{-1} \hat{S}_u .\]
 
{\bf 7.7. Proposition.} {\em The operators $\L_m+N\d_{m,0}/16$
satisfy the commutation relations $[\L_m,\L_n]$
$ = (m-n)\L_{m+n}$
of the algebra of vector fields on the line and annihilate the total
descendent potential of the semisimple Frobenius manifold: 
$\L_m\D =0,\ m=-1,0,1,2,...$. }  

\medskip

The formal Virasoro constraints described in the proposition
will be computed in the next section in the case of conformal 
Frobenius structures.

\bigskip
 
{\large \bf 8. Virasoro constraints.} 

\medskip

Consider a Frobenius structure of conformal dimension $D$ 
with the Euler vector field which is the sum of a linear diagonalizable and 
a constant vector field. 
Then in suitable flat coordinates of the metric it assumes the form
\[ E=\sum (1-d_{\a}) t^{\a}\p_{t^{\a}} + 
\sum_{\a:\ d_{\a}=1} \rho_{\a}\p_{t^{\a}}, \]
where the degree spectrum $\{ d_{\a} \}$ is confined on the interval
$[0,D]$ and is symmetric about $D/2$. The homogeneity conditions for
the calibration $S(1/z)$ and the asymptotical solution $T=\Psi R(z) e^{U/z}$
read respectively:
\[ (z\p_z + E ) S = \mu S + S (\mu+\rho/z)^* , \ \ 
   (z\p_z + E ) T = \mu T , \]
where $\mu=\operatorname{diag} (d_1-D/2,...,d_N-D/2)$
is anti-symmetric ($\mu = -\mu^* $), and $\rho = \rho^*$ is 
{\em $\mu$-nilpotent} in the sense of \cite{D1}. In 
Gromov -- Witten theory, $\mu $ is the 
{\em Hodge grading operator} and $\rho$
is the operator of multiplication by 
$c_1(TX)=\sum \rho_{\a} \p_{t^{\a}}$ in the {\em classical}
cohomology algebra.

\medskip

{\bf 8.1. Theorem.} {\em The Virasoro operators $7.7$ are 
$\L_m = \hat{L}_m^{\mu,\rho} + \frac{\d_{m,0}}{4} \operatorname{tr} (\mu\mu^*)$,
where $L_m^{\mu,\rho}$ are the infinitesimal symplectic transformations 
\[ L_m^{\mu,\rho} = z^{\mu}z^{-\rho} L_m z^{\rho}z^{-\mu} = 
z^{-1/2} (z\frac{d}{dz}z -\mu z +\rho)^{m+1} z^{-1/2} .\]} 

The Virasoro operators $8.1$ actually coincide with those introduced
in \cite{DZ1}. 

Regardless of the validity of the conjecture $7.6$,
our definition of the total descendent potential $\D$ reproduces 
correctly the formulas $5.4$ and $5.5+7.2$ for genus $0$ and $1$
descendent potentials in Gromov -- Witten theory. As a consequence,
we obtain a new proof of the main result in \cite{DZ1}:

\medskip

{\bf 8.2. Corollary.} {\em The genus $1$ Virasoro constraints hold true
for Gromov -- Witten invariants of almost K\"ahler manifolds with 
generically semisimple quantum cup-product.}

\medskip

{\bf 8.3.} {\em Remarks.} (a) Put $w=1/z$. Given a 
connection operator $D = d/dw + A(w)/w$ where
$A=A_0w+A_1w^2+...,\ A^*(-w)=-A(w)$, one can associate
to it a representation of the Lie algebra of vector
fields on the line to the algebra of quadratic
hamiltonians on $\H$ by taking 
$L^{D}_m:= w^{3/2} D^{m+1}w^{-1/2}$. Quantization
then yields a representation in the Fock space.  
If two connections are conjugated by a gauge transformation 
$D\mapsto S^{-1} D S $ where 
$S(w)=S_0+S_1w+...,\ S^*(-w)S(w)=\1$, then the 
corresponding representations are equivalent. 
According to \cite{D1}, the classes of gauge equivalence
(in the case when the residue operator $A_0, \ A_0=-A^*_0$,
is semisimple) are represented by the connections of the
form $d/dw + (\mu - w\rho )/w$ where $\rho $ is
polynomial in $w$ and is {\em $\mu$-nilpotent} 
(see \cite{D1}). The representation in the theorem corresponds to 
$D $ in the normal form with constant $\rho $.

(b) In the case of conformal Frobenius structures the
consistent PDE system $d \Phi = z^{-1} A(t) \Phi $ on $H$ 
can be completed, following \cite{D}, to a consistent 
PDE system  on $\CC \times H$ by adding the 
homogeneity equation 
$(d/dw + \mu/w -  E\bullet_t) \Phi =0 $.
The latter equation can be considered as a family $D_t$ of 
connections on $\CC $ depending on the parameter $t\in H$
and is {\em isomonodromic}. It is the role of the calibration 
$S_t(w)$ to conjugate the connection operators $D_t$ to the 
normal form $d/dw +\mu/w -\rho$ with the
fundamental solution $w^{-\mu}w^{\rho}$. 
In other words, $\Phi = S_t(w) w^{-\mu}w^{\rho}$ is
a fundamental solution to the PDE system on $\CC\times H$.
The key point in the proof of the theorem $8.1$ is that the 
asymptotical solution $T$ satisfies the same equations as
$\Phi $ and thus $S^{-1} T$ satisfies the same equations as 
$w^{-\mu}w^{\rho}$. 
  
\bigskip  

{\large \bf 9. Equivariant Gromov -- Witten theory.}
 
\medskip

Let the compact almost K\"ahler manifold $X$ be equipped with
a hamiltonian Killing action of a compact Lie group $G$. Then
the moduli spaces $X_{g,m,d}$ inherit the action, and the evaluation,
forgetting and contraction maps are $G$-equivariant. The construction
of the virtual fundamental class $[X_{g,m,d}]$ admits an equivariant
generalization \cite{R}. This allows one to introduce $G$-equivariant
counterparts of Gromov -- Witten invariants, gravitational descendents
and ancestors and of the generating functions $F^g_X, \F^g_X, \D_X, \A_X$,
etc. The invariants take values in $H^*_G(point; \QQ) = H^*(BG; \QQ)$,
the coefficient ring of equivariant cohomology theory. They reduce to
their non-equivariant versions under the restriction homomorphism
$H^*(BG; \QQ) \to H^*(point; \QQ)$.

The genus $0$ equivariant Gromov -- Witten potential 
$F^0_X$ defines on the equivariant cohomology space $H=H^*_G(X; \QQ [[Q]])$
the structure of a formal Frobenius manifold {\em over the ground ring}
$H^*_G(point; \QQ )$ (see \cite{Gi2}).  

Let us consider the case when $G$ is a torus $(S^1)^r$.
The cohomology algebra of the classifying space $BG=(\CC P^{\infty})^r$ 
is isomorphic to $\QQ [\l_1,...,\l_r]$. Here $\l_i \in Lie^*G$ are identified with 
a basis of infinitesimal characters of $1$-dimensional representations of 
$(S^1)^r$. The representations induce Hopf line bundles over 
$(\CC P^{\infty})^r$ whose Chern classes generate $H^*(BG)$. 

When the fixed points of the torus $G$ action on $X$ are isolated, 
the Frobenius structure on $H$ is semisimple. Indeed, even the classical 
equivariant cohomology algebra 
of $X$ with coefficients in the field of fractions $\QQ (\l )$ 
of the ground ring $H^*_G(point; \QQ)$ is semisimple. Therefore
the most of results of the previous sections apply in equivariant
Gromov -- Witten theory.  However, the Frobenius structure in question is 
not conformal: the Euler vector field expressing dimensional properties 
of Gromov -- Witten invariants is a derivation 
over $\QQ $ but, generally speaking, not over the ground ring $H^*_G(point; \QQ)$ 
since elements of this ring may have non-zero degrees. 
As a consequence, the part (c) of Proposition $6.7$ does not apply, \this
the construction $6.7 (a)$ of the series $R(z)$ is ambiguous, and the
ambiguity is as described in $6.7 (b)$.

Let $w_i,\ i=1,...,N$ be the fixed points of the torus $G$ action on $X$,
and $N_{l}^{(i)}$ denote the Newton polynomial $\sum_{j=1}^n \gx_j^{-l}(w_i)$
of the inverse infinitesimal characters $\gx_j^{-1}(w_i)$ of the torus $G$ 
action on the cotangent space $T^*_{w_i}X$ at the fixed point. We consider
$N_{l}^{(i)}$ as an element of the ground field $\QQ (\l_1,...,\l_r)$. 

\medskip

{\bf 9.1. Theorem.} {\em In the case of a torus $G$ action with
isolated fixed points, the total ancestor potential in the 
equivariant Gromov -- Witten theory on $X$ coincides with the total 
ancestor potential $6.8$ of the semisimple Frobenius structure on 
$H=H^*_G(X; \QQ [[Q]])$ (\this $e^{F^1_X}\A_X= e^C\hat{\Psi}\hat{R}\T$) provided that 
the series $R(z)$ is normalized by the condition that in the classical 
cohomology limit $Q\to 0$ it takes on the form $ R(z)|_{Q=0}= \exp \operatorname{diag} (b_1,...,b_N)$
where $b_i$ are the {\em Bernoulli series}
\[ b_i(z) = \sum_{k=1}^{\infty} N_{2k-1}^{(i)} \frac{B_{2k}}{2k}\frac{z^{2k-1}}{2k-1} .\] }

{\bf 9.2.} {\em Remarks.} (a) The theorem is equivalent (due to $5.1$) to the Theorem $2$
in \cite{Gi0} describing descendent potentials $\F^g_X$ with $g>1$ (plus a similar 
result of \cite{Gi1} about $F^1_X$). Both results are proved by fixed point localization 
in moduli spaces $X_{g,m,d}$ utilizing Kontsevich's technique of summations over graphs 
\cite{Kn1}. Justification of fixed point localization formulas for virtual fundamental
classes depends on the results of \cite{GrP} in the algebraic category and 
\cite{S1} in the general almost K\"ahler setting.
The Bernoulli series arise to offset the effect of Hodge integrals in localization formulas
by applying the theorem $4.1$ for $X=point$. 

(b) It is essential in the formulation $9.1$ that indices of canonical coordinates coincide
with the labels $i=1,...,N$ of fixed points. This is due to the 
phenomenon \cite{Gi2,Gi1} of {\em materialization} of canonical coordinates in fixed point 
localization theory, playing an essential role in the proofs \cite{Gi1, Gi0} otherwise too. 
The simplest instance of it is related to the equality 
\[ \sum  u^{i} = \ \text{the total "number" of elliptic curves in $X$ with given modulus} .\]
In the computation of the "number" via the fixed point technique, we can single out contributions
$v^{i}$ of those fixed curves where {\em the} elliptic irreducible component 
(with the required generic modulus) 
is mapped to the fixed point $w_i$. It turns out \cite{Gi2} that the equality $\sum v^{i} =
\sum u^{i}$ is {\em termwise} (which in particular establishes the aforementioned 
correspondence).          

(c) The ancestor potential $\A_X$ specializes to its non-equivariant counterpart in
the non-equivariant limit $\l =0$. It is not clear however how to derive from the theorem its non-equivariant 
version. It would suffice to show that the series $R(z)$ tends to its non-equivariant counterpart
(since all other ingredients of the formula obviously do). Yet, if for some
normalization of $R(z)$ the limit exists, then it does not exist for any other normalization.

\medskip

{\bf 9.3. Conjecture.} {\em When $X$ has generically semisimple non-equivariant quantum cohomology, 
the series $R(z)$ (normalized by the Bernoulli series as in $9.1$) has a non-equivariant limit at 
$\l = 0$ equal to the series $R$ from $7.6$.} 

\medskip

In the next section we prove this conjecture in the case of complex projective spaces.    
  
\bigskip

{\large \bf 10. Mirrors of complex projective spaces.} 

\medskip

According to the theorem $7.1$ the total descendent potential $\D$
of a semisimple Frobenius manifold is determined by the values of
$S, \Psi$ and $R$ at any semisimple point. We will exploit
the property of projective spaces and their products to
have semisimple {\em small} quantum cohomology algebra.

Let $G=T^{n+1}$ be the torus of diagonal unitary 
transformations of the standard Hermitian space 
$\CC^{n+1}$. We identify $H^*(BG,\QQ)$ with 
$\CC [\l_0,...,\l_n]$ and the equivariant cohomology
algebra $H^*_{G}(\CC P^n,\QQ)$ of the projectivized space
--- with $\CC [P,\l_0,...,\l_n]/((P-\l_0)...(P-\l_n)$.
Here $-P$ denote the equivariant $1$-st Chern class 
of the Hopf line bundle over $\CC P^n$. The multiplication
table in the basis $1, P,...,P^{n-1}$ in the {\em small}
equivariant quantum cohomology algebra of $\CC P^n$ has
the following well-known description:  
\[ (P-\l_0)\bullet ... \bullet (P-\l_n) = Q,\ \ 
\text{while}\ \ P\bullet P^{k} = P^{k+1} \ \text{for} \ 
k<n. \] 
Here the generator $Q$ in the Novikov ring is identified with $\exp t$ 
where $t$ is the coordinate on the parameter space $H^2(\CC P^n, \CC)$ 
of the {\em small} quantum cohomology algebra. (The identification is
possible due to the {\em divisor equation}, see for instance \cite{Gi2}.) 
Respectively, the connection $\nabla_z $ from $6.7$
restricted to this parameter space yields
the following system of linear differential equations:
\[ z \frac{d}{dt} I_k = I_{k+1}\ 
\text{for $k<n$, and}\ \ \ 
(z \frac{d}{dt}-\l_0)...(z \frac{d}{dt}-\l_n)\ I_0 \ 
=\ e^t\ I_0.\]

Introduce the complex oscillating integral
\[ I = \int_{\Gamma \subset \{ x_0...x_n=Q \}}  
e^{(x_0+...+x_n)/z}\ x_0^{\l_0/z}...\ x_n^{\l_n/z}
\ \frac{dx_0\w ...\w dx_n}{d(x_0...dx_n)} .\]
The cycles $\Gamma$ in the oscillating integrals
of the form $\int_{\Gamma} e^{F(x)/z}\phi (x) dx$
may be non-compact and are constructed by means of Morse 
theory for the functions $\operatorname{Re} (F/z)$.

\medskip

{\bf 10.1. Theorem.} \ \ $ (z Q\frac{d}{dQ}-\l_0)...
(z Q\frac{d}{dQ} -\l_n)\ I \ =\ Q\ I $.

\medskip

{\em Proof.} It is convenient to rewrite the integral in
logarithmic coordinates: $Q=e^t, \ x_i=e^{T_i}$:
\[ I=\int_{ \Gamma \subset \{ t=\sum T_i \} }
e^{\sum (e^{T_i}+\l_i T_i) /z} \ \frac{dT_0\w ... \w dT_n}
{d(T_0+...+T_n)} .\]
The projection $(T_0,...,T_n)\mapsto t=T_0+...+T_n$ maps
each $\p_{T_i}$ to $\p_t$. Using this we conclude that
application of $z d/dt -\l_i$ to the integral can be
replaced by multiplication of the integrand by $e^{T_i}$.
Doing this consecutively for $i=0,1,...,n$ we find 
the integrand multiplied by $e^{T_0}...e^{T_n}=e^t$.
Thus $(z d/dt-\l_0)...(zd/dt-\l_n) I = e^t I$.

\medskip

{\bf 10.2.} {\em Remark.} The result is an equivariant
version of the mirror theorem for complex projective spaces 
we reported in Summer $93$ at seminars in Lyon, 
Strasbourg and Oberwolfach. Many mathematical results inspired by or 
predicted on the basis of mirror theory have
been proved since then. Yet the argument presently discussed 
seems to be one of the first formal 
mathematical applications of a mirror theorem.

\medskip

An asymptotical solution 
$T=\Psi R(z) \exp (U/z)$ to the system $\nabla_z T=0$
can be obtained from stationary
phase asymptotics near the critical points of the phase function
of the integral $I$ and its derivatives.

\medskip

{\bf 10.3.} {\em Example.} We illustrate the procedure of expanding a
complex oscillating integral near a non-degenerate critical point of the
phase function with the critical value $u$: 
\[ \int e^{(u-x^2+\a x^3 + ...)/z }(\b + \c x + ...) dx = 
\sqrt{z} e^{u/z} \int e^{-y^2} e^{\a y \sqrt{z}+...}(\b + \c y \sqrt{z}+...) dy .\]
Discarding the dimensional factor $\sqrt{z}$ in front of the integral and computing 
momenta of the Gaussian distribution yields an asymptotical expansion of the form 
$ e^{u/z} \  (a+bz+cz^2+...) $.

\medskip

In the oscillating integral $I$, the critical points of the phase function 
are constrained extrema of the function 
\[ e^{T_0}+...+e^{T_n}+\l_0T_0+...+\l_nT_n-P(T_0+...+T_n-t) \]
with the Lagrange multiplier $P$.  They satisfy $e^{T_{\a}} =P-\l_i$ where 
$(P-\l_0)...(P-\l_n)=Q$. 
Near $Q=0$ and generic $\l$, the $n+1$ roots $P^{(i)}(Q),
\ i=0,...,n$, to the latter equation are distinguished 
by their values $P^{(i)}(0)=\l_i$. 
On the other hand the basis $\{ P^k, k=0,...,n \}$ in $H^*(\CC P^n)$
corresponds to the basis of integrals $I_k = z (d/dt)^k I$ (since $P^{\bullet k}=P^k$
for $k\leq n$). In the asymptotical solution
$\Psi R(z) \exp (U/z)$ the diagonal matrix $U$ consists of the critical 
values $u^i$, and the matrix entries of $\Psi R(z)$ are obtained from 
asymptotical expansions of $I_k$ near the critical point $P^{(i)}$.

The asymptotical solution thus constructed 
automatically admits the non - equivariant limit $\l = 0$ when $Q\neq 0$
since the phase functions, amplitudes and the critical points in the oscillating integrals
depend continuously on $\l $. For generic $\l $, the solution depends continuously on $Q$ 
up to $Q=0$. 

At $Q=0$ the matrix $\Psi $ describes the transition from the basis $\{ P^k \}$ to
the basis    
$L_j(P)=\prod_{\a} (P-\l_{\a})/(\l_j-\l_{\a})^{1/2}$ of (suitably normalized) Lagrange
interpolation polynomials. Indeed, the Lagrange interpolation polynomials represent 
the basis of delta-functions of
of the fixed points $w_j, j=0,...,n$, in the equivariant cohomology algebra of $\CC P^n$   
and therefore coincide with the canonical idempotents of the semisimple algebra $(T_uH,\bullet_u)$ 
in the classical cohomology limit $Q=0$.        
Thus the matrix elements of the series $R(z)$ at $Q=0$ are extracted from
asymptotical expansions near the critical points $P^{(i)}$ of the integrals 
$L_j (z d/dt) I$ in this limit.

For symmetry reasons it suffices to consider only the 
critical point corresponding to $P^{(0)}$. Computing $L_j(z d/dt) I$ 
as in the proof of $10.1$ 
we find an extra factor proportional to $Q/x_j$ in the integrand of the oscillating integral $I$.
We rewrite the integral in the chart $(x_1,...,x_n)$
where $x_0=Q/x_1...x_n$. It still contains $P^{(0)}$ 
in the limit $Q=0$ when $x_0 = P^{(0)}-\l_0$ vanishes.  
We find $L_j(zd/dt) I$ proportional to 
\[ e^{\l_0\ln Q/z} \int e^{[\frac{Q}{x_1...x_n} +x_1+...+x_n 
+(\l_1-\l_0)\ln x_1 +...+(\l_n-\l_1)\ln x_n]/z }\ 
\frac{Q}{x_j}\ \frac{dx_1\w ...\w dx_n}{x_1...x_n} .\]
We see that the integral vanishes at $Q=0$ unless $j=0$. This agrees
with the fact of our general theory that $R$ is diagonal
at $Q=0$. 
For $j=0$ we have $Q/x_0=x_1...x_n$, and the integral at 
$Q=0$ factors into $1$-dimensional ones:
\[ \int e^{(\sum_{i=1}^n x_i+(\l_i-\l_0)\ln x_i)/z} \ dx_1\w...\w dx_n = 
\prod_{i=1}^n \int e^{x/z} x^{(\l_i-\l_0)/z} dx .\]
The stationary phase asymptotics of the latter integral is
the same as for the product 
$\Gamma (1+(\l_1-\l_0)/z) ... \Gamma (1+(\l_n-\l_0)/z)$
of the gamma-functions as $z\to 0$ since  
\[ \int e^{x/z}\ x^{\l /z}\ dx =
(-z)^{1+\l/z} \int e^{-v}\ v^{\l/z}\ dv \ 
\sim \ \Gamma (1+\frac{\l}{z}) .\]

\medskip

{\bf 10.4. Lemma.} \ \ $ \ln \Gamma (1+s) \sim s\ln (s/e) +
\frac{1}{2}\ln (2\pi s) + 
\sum_{k=1}^{\infty} \frac{B_{2k}}{2k} 
\frac{s^{1-2k}}{2k-1}$.   

\medskip

{\em Proof.} This is well-known, see for instance 
\cite{JEL}.   

\medskip 

Subtracting the critical value (Stirling's 
approximation) we find that at $Q=0$ 
\[ \ln R^{00}=\sum_{k=1}^{\infty} \frac{B_{2k}}{2k}
\frac{z^{2k-1}}{2k-1} \sum_{\a\neq 0} 
(\l_{\a}-\l_0)^{1-2k} .\]
This coincides with $ b_0 (z)$ in $9.1$ since $\l_{\a}-\l_0$ are
exactly the infinitesimal characters of the torus $G$
action on $T^*_{w_0}\CC P^n$. 

\medskip

Finally, our computation carries over without significant changes
to arbitrary Fano toric manifolds $X$. It is essential here that 
(i) Fano toric manifolds have semisimple {\em small} quantum cohomology
(as it follows from the explicit description of such cohomology
given by V. Batyrev \cite{Ba}) and (ii) {\em equivariant} genus $0$ 
GW-invariants of Fano toric manifolds admit a mirror 
description \cite{Gi1} by complex oscillating integrals generalizing
Theorem $10.1$. 

Combining the computation with $5.1$ and $9.1$ 
(or equivalently --- with \cite{Gi0}) we prove the conjecture $7.6$.
 
\medskip  

{\bf 10.5. Theorem.} {\em The total descendent potential of a toric
Fano manifold $X$ is given, up to a non-zero
constant factor, by the formula
$ \D_X = e^{C(t)}\hat{S_t} \hat{\Psi_t} \hat{R_t} \T$ 
where $t\in H^2(X;\CC)$.} 

\medskip

{\bf 10.6. Corollary.} {\em The total descendent potential $\D_X$ 
of a toric Fano manifold satisfies $(\L_m +N\d_{m,0}/16) \D_X =0$ for
$m=-1,0,1,2,...$, where $\L_m$ are the Virasoro operators $8.1$.}

\medskip

{\bf 10.7.} {\em Remarks.} 
(a) The corollary was conjectured by T. Eguchi, 
K. Hori, M. Jinzenji, C.-S. Xiong  and S. Katz \cite{EHX, EJX}. 

(b) The so called {\em $\l_g$-conjecture} and some other non-trivial properties 
of Hodge integrals over Deligne -- Mumford spaces
were discovered by Getzler and Pandharipande \cite{GeP} on the basis of Virasoro conjecture for $\CC P^1$. The Corollary therefore backs 
up the original arguments of \cite{GeP} which appear to be
more natural than the proof of the properties found 
later by Faber and Pandharipande \cite{FP}.

(c) The role
the gamma-function $10.4$ plays in a general theorem 
$9.1$ should probably remind us  
about the relationship of quantum cohomology 
theory with its heuristic roots \cite{Gi4,Gi3} ---
$S^1$-equivariant Floer theory on loop spaces.
In fact the gamma-function $10.4$ will make a more
systematic appearance in \cite{CG} in the context of general
``quantum'' versions of the Riemann -- Roch formula,
Serre duality and Lefschetz hyperplane section theorem. 

\medskip

{\large \bf  Acknowledgments.} The author is thankful to 
H. Chang, T. Coates, Y.-P. Lee and R. Pandharipande
for numerous stimulating discussions and to J. Morava for 
pointing at his notes \cite{Mo} which indicate how 
the formalism of Fock spaces discussed above emerges naturally 
from $S^1$-equivariant cobordism theory of loop spaces.

\enddocument